\documentclass[12pt]{amsart}
\usepackage{amsmath}
\usepackage{amsfonts}
\usepackage{amssymb}
\usepackage{epsfig,latexsym,graphicx}

\newtheorem{Thm}{Theorem}[section]
\newtheorem{thm}[Thm]{Theorem}

\newtheorem{lemma}[Thm]{Lemma}

\newtheorem{remark}[Thm]{Remark}

\newtheorem{cor}[Thm]{Corollary}

\newcommand{\ip}[2]{{\langle#1,#2\rangle}}
\newcommand{\norm}[1]{{\|#1\|}}

\def\ldots{\mathinner{\ldotp\ldotp\ldotp}}
\def\ldots{\mathinner{\cdotp\cdotp\cdotp}}

\def \K{\rm {\bf K}}

\def \cal{\mathbb}

\def \beq{\begin{eqnarray*}}
\def \eeq{\end{eqnarray*}}

\def \R{\mathbb{R}}
\def \C{\mathbb{C}}
\def \K{\mathbb{K}}
\def \E{\mathbb{\cal E}}

\newcommand{\fc}{{\mathbb F}}
\newcommand{\I}{{\mathbb I}}
\newcommand{\deltavect}{v}

\title{Reconstruction of Signals from Magnitudes of Redundant Representations}
\author[R.~Balan]{Radu Balan}
\address{Department of Mathematics \\ University of Maryland, College Park MD 20742}
\email[R.~Balan]{rvbalan@math.umd.edu} 

\begin{document}

\begin{abstract}
This paper is concerned with the question of reconstructing
a vector in a finite-dimensional real or complex Hilbert space when 
only the magnitudes of the coefficients of the vector under a redundant linear map
are known. We present new invertibility results as well an iterative algorithm that finds the least-square solution
and is robust in the presence of noise. We analyze its numerical performance by
comparing it to two versions of the Cramer-Rao lower bound.
\end{abstract}

\maketitle

\section{Introduction}

This paper is concerned with the question of reconstructing
a vector $x$ in a finite-dimensional real or complex Hilbert space $H$ when 
only the magnitudes of the coefficients of the vector under a redundant linear map
are known.

Specifically our problem is to reconstruct $x\in H$ up to a global phase factor from the magnitudes 
$\{ |\ip{x}{f_k}|~,~1\leq k\leq m\}$ where $\{f_1,\ldots,f_m\}$ is a frame (complete system) for $H$.

A previous paper \cite{BCE06} described the importance of this problem
to signal processing, in particular to the analysis of speech. 
Of particular interest is the case when the coefficients are obtained
from a Windowed Fourier Transform (also known as Short-Time Fourier
Transform), or an Undecimated Wavelet Transform (in audio and image
signal processing).  While \cite{BCE06} presents some necessary and
sufficient conditions for reconstruction, the general problem of finding
fast/efficient algorithms is still open.  In \cite{Bal2010} we describe one
solution in the case of STFT coefficients.

For vectors in real Hilbert spaces,  the reconstruction problem 
is easily shown to be equivalent to a
combinatorial problem.  In \cite{BCE07} this problem is
further proved to be equivalent to a (nonconvex) optimization problem.

A different approach (which we called the "algebraic approach") was proposed
in \cite{Bal2009}. While it applies to both real and complex cases, noisless
and noisy cases, the approach requires solving a linear system of size exponentially 
in space dimension.  The algebraic approach mentioned earlier
generalizes the approach in \cite{BBCE07} where reconstruction is performed
with complexity $O(n^2)$ (plus computation of the principal eigenvector for a
matrix of size $n$). However this method requires $m=O(n^2)$ frame vectors.

Recently the authors of \cite{CSV12} developed a convex
optimization algorithm (PhaseLift) and proved its ability to perform exact reconstruction 
in the absence of noise, as well as its stablity under noise conditions. 
In a separate paper, \cite{CESV12}, the authors developed further a similar algorithm in the case
of windowed DFT transforms. 

In this paper we analyze an iterative algorithm based on regularized least-square criterion.
The organization of the paper is as follows. In section \ref{sec2} we define the problem
explicitely. In section \ref{sec3} we describe our approach and prove some convergence results.
In section \ref{sec4} we establish the Cramer-Rao lower bound for benchmarking its
performance which is analyzed in section \ref{sec5}. Section \ref{sec6} contains
conclusions and is followed by references.

\section{Background\label{sec2}}

Let us denote by $H$ the n-dimensional Hilbert space $\R^n$ (the real case) or $\C^n$
 (the complex case) with scalar product $\ip{}{}$. Let $\fc=\{f_1,\ldots,f_m\}$ be a spanning set of $m$ vectors in $H$.
In finite dimension (as it is the case here) such a set forms a {\em frame}. In the infinite
dimensional case, the concept of frame involves a stronger property than completness (see 
for instance \cite{Cass-artofframe}). 
We review additional terminology and properties which remain still true in the infinite dimensional setting.
The set $\fc$ is frame if and only if there are two positive constants $0<A\leq B<\infty$ (called frame bounds)
so that
\[ A\norm{x}^2 \leq \sum_{k=1}^m |\ip{x}{f_k}|^2 \leq B\norm{x}^2 \]
When we can choose $A=B$ the frame is said {\em tight}. For $A=B=1$ the frame is called {\em Parseval}.
A set of vectors $\fc$ of the $n$-dimensional Hilbert space $H$ is said to have {\em full spark} if any subset
of $n$ vectors is linearly independent. 

For a vector $x\in H$, the collection of coefficients
$\{ \ip{x}{f_k}~,~1\leq k\leq m\}$ represents the analysis of vector $x$ given by the frame $\fc$.
In $H$ we consider the following equivalence relation:
\begin{equation}
\label{eq:equiv}
x,y\in H ~~,~~ x\sim y~~{\rm iff}~~y=zx~{\rm for~some~scalar~z~with}~|z|=1
\end{equation}
Let $\hat{H}=H/\sim$ be the set of classes of equivalence induced by this relation.
Thus $\hat{x}=\{x,-x\}$ in the real case (when $H=\R^n$), and $\hat{x}=\{e^{i\alpha}x,
0\leq\alpha <2\pi \}$ in the complex case (when $H=\C^n$). The analysis map induces
the following nonlinear map
\begin{equation}
\label{eq:varphi}
\varphi :\hat{H}\rightarrow (\R^{+})^m~~,~~\varphi(\hat{x})=(|\ip{x}{f_k}|^2)_{1\leq k\leq m}
\end{equation}
where $\R^{+}=\{x~~,~~x\in\R~,~x\geq 0\}$ is the set of nonnegative real numbers.
In previous papers \cite{BCE06} we studied when the nonlinear map $\varphi$ is injective. We review these results below.
In this paper we describe an algorithm to solve the equation
\begin{equation}
\label{eq:eq}
y = \varphi(x)
\end{equation}
and then we study its performance in the presence of additive white Gaussian noise when the model becomes
\begin{equation}
\label{eq:model}
y = \varphi(x) + \nu ~~,~~\nu\sim {\cal N}(0,\sigma^2)
\end{equation}
We shall derive the Cramer-Rao Lower Bound (CRLB) for this model and compare its performance to this bound.

\subsection{Existing Results}

We revise now existing results on injectivity of the nonlinear map $\varphi$. In general a subset $Z$ of a topological space 
is said {\em generic} if its open interior is dense. However in the following statements, the term
{\em generic} refers to Zarisky topology: a set $Z\subset \K^{n\times m}= \K^{n}\times\cdots\times \K^n$ is said {\em generic}
if $Z$ is dense in $\K^{n\times m}$ and its complement is a finite union of zero sets of polynomials in $nm$ variables 
with coefficients in the field $\K$ (here $\K=\R$ or $\K=\C$). 
\begin{thm}[\cite{BCE06}Th.2.8] \label{th2.1} In the real case when $H=\R^n$ the following are equivalent:
\begin{enumerate}
\item The nonlinear map $\varphi$ is injective;
\item For any disjoint partition of the frame set $\fc=\fc_1\cup\fc_2$, either $\fc_1$ spans $H$ or $\fc_2$ spans $H$.
\end{enumerate}
\end{thm}
\begin{cor}[\cite{Fink04}Th.I;\cite{BCE06}Th.2.2,Prop.2.5,Cor.2.6] \label{th2.2} The following hold true in the real case $H=\R^n$:
\begin{enumerate}
\item If $\varphi$ is injective then  $m\geq 2n-1$;
\item If $m\leq 2n-2$ then $\varphi$ cannot be injective;
\item If $m=2n-1$ then $\varphi$ is injective if an only if $\fc$ is full spark;
\item If $m\geq 2n-1$ and $\fc$ is full spark then the map $\varphi$ is injective;
\item If $m\geq 2n-1$ then for a generic frame $\fc$ the map $\varphi$ is injective.
\end{enumerate}
\end{cor}
\begin{thm}[\cite{Fink04}Th.II;\cite{BCE06}Th.3.3] \label{th2.3} In the complex case when $H=\C^n$ the following hold true:
\begin{enumerate}
\item If $m\geq 4n-2$ then for a generic frame $\fc$ the map $\varphi$ is injective;
\item If $\varphi$ is injective then $m\geq 3n-2$;
\item If $m\leq 3n-3$ then the map $\varphi$ cannot be injective.
\end{enumerate}
\end{thm}

\subsection{New Injectivity Results}

We obtain equivalent conditions to (2) in Theorem \ref{th2.1}. In the real case these new conditions are
equivalent to $\varphi$ being injective; in the complex case they are only necessary condition for injectivity.

\begin{thm}\label{th2.4} Given a $m$-set of vectors $\fc=\{f_1,\ldots,f_m\}\subset H$ the following conditions are
equivalent:
\begin{enumerate}
\item  For any disjoint partition of the frame set $\fc=\fc_1\cup\fc_2$, either $\fc_1$ spans $H$ or $\fc_2$ spans $H$;
\item For any two vectors $x,y\in H$ if $n\neq 0$ and $y\neq 0$ then 
$$\sum_{k=1}^m |\ip{x}{f_k}|^2|\ip{y}{f_k}|^2 >0;$$
\item There is a positive real constant $a_0>0$ so that for all $x,y\in H$,
\begin{equation}
\label{eq:q4}
\sum_{k=1}^m |\ip{x}{f_k}|^2 |\ip{y}{f_k}|^2 \geq a_0 \norm{x}^2 \norm{y}^2
\end{equation}
\item There is a positive real constant $a_0>0$ so that for all $x\in H$,
\begin{equation}
\label{eq:q2}
R(x) := \sum_{k=1}^m |\ip{x}{f_k}|^2 \ip{\cdot}{f_k}f_k \geq a_0 I
\end{equation}
where the inequality is in the sense of quadratic forms.
\end{enumerate}
\end{thm}
\begin{remark}
The constants in (3) and (4) above are the same (hence the same notation).
\end{remark}

\emph{Proof}

$(1)\Rightarrow(2)$. We prove by contradiction: $\neg (2)\Rightarrow\neg(1)$. Assume there are $x,y\in H$, $x\neq 0$, $y\neq 0$
so that $\sum_{k=1}^m |\ip{x}{f_k}|^2|\ip{y}{f_k}|^2=0$. Then $\ip{x}{f_k} \ip{y}{f_k}=0$ for all $k$. Let $I_{x}=\{k~,~
\ip{x}{f_k}=0\}$, and set $\fc_1=\{f_k~,~k\in I_x\}$. Let $I_y=\{1,\ldots,m\}\setminus I_x$, and set $\fc_2=\fc\setminus\fc_1$.
Since $x$ is orthogonal to $\fc_1$ it follows that $\fc_1$ cannot span the whole $H$; similarly $\fc_2$ cannot span $H$ because
$y$ is orthogonal to all $\fc_2$. This contradicts $(1)$.

$(2)\Rightarrow(3)$. The unit sphere $S_1(H)$ is compact in $H$ and so is $S_1(H)\times S_1(H)$. Since the map
\[ (x,y) \mapsto \sum_{k=1}^m |\ip{x}{f_k}|^2 |\ip{y}{f_k}|^2 \]
is continuous, it follows
\begin{equation}
\label{eq:a0}
 a_0 := min_{(x,y)\in S_1(H)\times S_1(H)} \sum_{k=1}^m |\ip{x}{f_k}|^2 |\ip{y}{f_k}|^2 > 0. 
\end{equation}
By homogeneity for any $x,y\in H$, $x\neq 0,y\neq 0$ we obtain:
\[ \sum_{k=1}^m |\ip{x}{f_k}|^2|\ip{y}{f_k}|^2 = \norm{x}^2\norm{y}^2 \sum_{k=1}^m |\ip{\frac{x}{\norm{x}}}{f_k}|^2
|\ip{\frac{y}{\norm{y}}}{f_k}|^2 \geq a_0\norm{x}^2\norm{y}^2 \]
hence (\ref{eq:q2}). If either $x=0$ or $y=0$ then (\ref{eq:q2}) holds true.

$(3)\Rightarrow(4)$. Follows immediately by definition of quadratic forms.

$(4)\Rightarrow(1)$. We prove by contradiction: $\neg(1)\Rightarrow\neg(4)$. If there is a partition $\fc=\fc_1\cup\fc_2$
so that neither $\fc_1$ spans $H$ nor $\fc_2$ spans $H$, then there are two non-zero vectors $x,y\in H$ so that $x\perp\fc_1$
and $y\perp\fc_2$. Thus $\ip{x}{f_k}\ip{y}{f_k}=0$ for all $k$. In turn this means $y\in ker(R(x))$ which contradicts (\ref{eq:q2}).
~~$\Box$.

Note the proof of this result produced the following condition equivalent to negating any of the statements of Theorem \ref{th2.4}:
There are two non-zero vectors $x,y\in H$ and a subset $\fc'\subset\fc$ so that $\ip{x}{f}=0$ for all $f\in\fc'$, and
$\ip{y}{f}=0$ for all $f\in\fc\setminus\fc'$. Then one can immediatly check that $x+y$ and $x-y$ are two non-equivalent
vectors in $H$ with respect to the (real or complex) equivalence relation $\sim$, and yet $\varphi(x+y)=\varphi(x-y)$; hence $\varphi$ cannot
be injective. We thus obtained
\begin{cor}\label{th2.5}
\begin{enumerate}
\item When $H=\R^n$, $\varphi$ is injective if and only if any (and hence all) of the conditions of Theorem \ref{th2.4} is satisfied;
\item When $H=\C^n$ , if $\varphi$ is injective then all conditions of Theorem \ref{th2.4} must hold.
\end{enumerate}
\end{cor}

\section{Our approach: Regularized Iterative Least-Square Optimization\label{sec3}}

Consider the additive noise model in (\ref{eq:model}). Our data is the vector $y\in\R^m$. Our goal is to find
a $x\in H$ that minimizes $\norm{y-\varphi(x)}$, where we use the Euclidian norm. As discussed also in section \ref{sec4},
the least-square error minimizer represents the Maximum Likelihood Estimator (MLE) when the noise is Gaussian.
In this section we discuss an optimization algorithm for this criterion.
Consider the following function
\begin{eqnarray}\label{eq:J}
&& J:H\times H\times \R^{+}\times\R^{+}\rightarrow \R^{+} \\
J(u,v,\lambda,\mu) & = & \sum_{k=1}^m|y_k - \ip{u}{f_k}\ip{f_k}{v}|^2+
\lambda\norm{u}^2 +\mu\norm{u-v}^2+\lambda\norm{v}^2. \nonumber
\end{eqnarray}
Our goal is to minimize $\norm{y-\varphi(u)}^2=J(u,u,0,\mu)$ over $u$, for some (and hence any) value $\mu\in\R^{+}$.
Our strategy is the following iterative process:
\begin{equation}\label{eq:xt1}
x^{t+1} = argmin_u J(u,x^{t},\lambda_t,\mu_t)
\end{equation}
for some initialization $x^{0}$ and policy $(\lambda_t)_{t\geq 0}$ and $(\mu_t)_{t\geq 0}$.

\subsection{Initialization}
Consider the regularized least-square problem:
\[ min_u J(u,u,\lambda,0) = min_u \norm{y-\varphi(u)}^2 + 2\lambda\norm{u}^2 \]
Note the following relation
\begin{eqnarray}
J(u,u,\lambda,0) &=& \norm{y}^2 +2\lambda\norm{u}^2 - 2\sum_{k=1}^m y_k |\ip{u}{f_k}|^2 + \sum_{k=1}^m
|\ip{u}{f_k}|^4 \nonumber \\
& = & \norm{y}^2 + 2\ip{(\lambda I-Q)u}{u}+\sum_{k=1}^m |\ip{u}{f_k}|^4 \label{eq:J0}
\end{eqnarray}
where
\begin{equation}
\label{eq:Q}
 Q=\sum_{k=1}^m y_k \ip{\cdot}{f_k}f_k. 
\end{equation}
For $\lambda>\norm{Q}$ the optimal solution is $u=0$. 
Note that if  $Q\leq 0$ as a quadratic form then the optimal solution of $min_u \norm{y-\varphi(u)}^2$ is $u=0$.
Consequently we assume the largest eigenvalue of $Q$ is positive. As $\lambda$ decreases
the optimizer remains small. Hence we can neglect the forth order term in $u$ in the expansion above and obtain:
\[ J(u,u,\lambda,0) \approx \norm{y}^2 + 2\ip{(\lambda I-Q)u}{u} \] 
Thus the critical value of $\lambda$ for which we may get a nonzero solution is $\lambda=max eig(Q)$ is the
maximum eigenvalue of $Q$. Let us denote by $e_1$ this (positive) eigenvalue and $v_1$ its associated normalized eigenvector.
This suggests to initialize $\lambda=\alpha e_1$ for some $\alpha\leq 1$ and $x^{0}=\beta v_1$, for some scalar $\beta$.
Substituting into (\ref{eq:J0}) we obtain
\[ J(\beta v_1,\beta v_1,\alpha e_1,0) = \norm{y}^2 - 2(1-\alpha)e_1\beta^2 + (\sum_{k=1}^m |\ip{v_1}{f_k}|^4)\beta^4 \]
For fixed $\alpha$, the minimum over $\beta$ is achieved at
\begin{equation}\label{eq:beta0}
\beta_0=\sqrt{\frac{(1-\alpha)e_1}{\sum_{k=1}^m|\ip{v_1}{f_k}|^4}}
\end{equation}
The parameter $\mu$ controls the step size at each iteration. The larger the value the smaller the step. On the other
hand, a small value of this parameter may produce an unstable behavior of the iterates. In our implementation we use
the same initial value for both $\lambda$ and $\mu$:
\begin{equation}
\label{eq:mu0}
\mu_0 = \lambda_0 = \alpha e_1
\end{equation}
 
\subsection{Iterations}
Optimization problem (\ref{eq:xt1}) admits a closed form solution. Specifically, expanding the quadratic in $u$ we obtain
\begin{eqnarray*}
 J(u,x^t,\lambda_t,\mu_t) & = & \ip{(R_t+\lambda_t+\mu_t)u}{u} - \ip{u}{(Q+\mu_t)x^t} - \ip{(Q+\mu_t)x^t}{u} + \\
 & + & \norm{y}^2 + (\lambda_t+\mu_t)\norm{x^t}^2 
\end{eqnarray*}
where
\begin{equation}\label{eq:Rt}
R_t = \sum_{k=1}^m |\ip{x^t}{f_k}|^2 \ip{\cdot}{f_k}f_k
\end{equation}
and $Q$ is defined in (\ref{eq:Q}). We obtain that $x^{t+1}$ satisfies the following linear equation
\begin{equation}
\label{eq:xt+1}
(R_t+\lambda_t+\mu_t)x^{t+1} = (Q+\mu_t)x^t
\end{equation}
Note the quadratic form on the left hand side is bounded below by
\[ R_t+\lambda_t+\mu_t \geq a_0\norm{x^t}^2 + \lambda_t+\mu_t > 0 \]
where $a_0$ is given by (\ref{eq:a0}).

\subsection{Convergence}
Denote $j_t=J(x^{t+1},x^t,\lambda_t,\mu_t)=min_u J(u,x^t,\lambda_t,\mu_t)$. We have the following general result:
\begin{thm}\label{th3.1}
Assume $\lambda_0\geq\lambda_1\geq\cdots\geq\lambda_t\geq\cdots$ and $\mu_0\geq\mu_1\geq\cdots\geq\mu_t\geq\cdots$.
Then for any initialization $x^0$ the sequence $(j_t)_{t\geq 0}$ is monotonically decreasing and therefore convergent.
\end{thm}
This theorem follows immediately from the following lemma:
\begin{lemma}\label{lem3.2}
Assume $\lambda_t\geq\lambda_{t+1}$ and $\mu_t\geq\mu_{t+1}$, Then $j_t\geq j_{t+1}$.
\end{lemma}

{\em Proof} First remark the symmetry
\begin{equation}
J(u,v,\lambda,\mu)=J(v,u,\lambda,\mu).
\end{equation}
Then we have:
\begin{eqnarray*}
 j_{t+1} & = & J(x^{t+2},x^{t+1},\lambda_{t+1},\mu_{t+1}) \leq J(x^{t},x^{t+1},\lambda_{t+1},\mu_{t+1}) \\
 & \leq & J(x^t,x^{t+1},\lambda_t,\mu_t) = J(x^{t+1},x^t,\lambda_t,\mu_t) = j_t 
\end{eqnarray*}
This concludes the proof of the lemma. $\Box$

\subsection{First Algorithm}\label{subsec3.4}
We are now ready to state the first optimization algorithm:

{\bf Input data}: Frame set $\fc=\{f_k~,~1\leq k\leq m\}$, measurements $y=(y_k)_{k=1}^m$, initialization
parameter $\alpha$, stopping criterion $\epsilon$, or maximum number of iterations $T_{max}$.

{\bf Initialization}: Compute matrix $Q$ in (\ref{eq:Q}) and its principal eigenvalue $e_1$ and eigenvector $v_1$. 
Compute $\beta_0$ in (\ref{eq:beta0}). Set $t=0$ and
\[ \lambda_0=\mu_0=\alpha e_1~~,~~x^0 =\beta_0 v_1 \]

{\bf Iterate}. Repeat:
\begin{enumerate}
\item Compute $R_t$ given by (\ref{eq:Rt});
\item Solve (\ref{eq:xt+1}) for $x^{t+1}$;
\item Update $\lambda_t$, $\mu_t$ using a preset or adaptive policy (more details are provided in section \ref{sec5});
\item Compute $j_t=J(x^{t+1},x^t,\lambda_t,\mu_t)$ and increment $t=t+1$;
\end{enumerate}
Until $t>T_{max}$, or $j_{t-1}-j_t < \epsilon$, or $\lambda_t<\epsilon$.

{\bf Outputs}: Estimated signal $x^{t}$, criterion $j_t$, error $\norm{y-\varphi(x^t)}^2$.

\subsection{Second Algorithm\label{subsec3.5}}
Results of numerical simulations (see section \ref{sec5}) suggest the adaptation of $\lambda$ and $\mu$ is too agressive.
Instead of running the algorithm until $\lambda< 1e-8$ (a small value), we implemented a second algorithm where we track
the mean-square error:
\[ L(x^t) = \sum_{k=1}^m |y_k - |\ip{x^t}{f_k}|^2|^2 \]
and return the minimum value. We thus obtain a second algorithm:

{\bf Input data}: Frame set $\fc=\{f_k~,~1\leq k\leq m\}$, measurements $y=(y_k)_{k=1}^m$, initialization
parameter $\alpha$, stopping criterion $\epsilon$, or maximum number of iterations $T_{max}$.

{\bf Initialization}: Compute matrix $Q$ in (\ref{eq:Q}) and its principal eigenvalue $e_1$ and eigenvector $v_1$. 
Compute $\beta_0$ in (\ref{eq:beta0}). Set $t=0$ and
\[ \lambda_0=\mu_0=\alpha e_1~,~x^0 =\beta_0 v_1~,~L_{optim} = \sum_{k=1}^m |y_k - |\ip{x^0}{f_k}|^2|^2~,~x_{optim}=x^0 \]

{\bf Iterate}. Repeat:
\begin{enumerate}
\item Compute $R_t$ given by (\ref{eq:Rt});
\item Solve (\ref{eq:xt+1}) for $x^{t+1}$;
\item Update $\lambda_t$, $\mu_t$ using a preset or adaptive policy (more details are provided in section \ref{sec5});
\item Compute $j_t=J(x^{t+1},x^t,\lambda_t,\mu_t)$; 
\item Compute $L_{t+1}=\sum_{k=1}^m |y_k - |\ip{x^{t+1}}{f_k}|^2|^2$;
\item If $L_{t+1}<L_{optim}$ then $L_{optim}=L_{t+1}$ and $x_{optim}=x^{t+1}$;
\item increment $t=t+1$;
\end{enumerate}
Until $t>T_{max}$, or $j_{t-1}-j_t < \epsilon$, or $\lambda_t<\epsilon$.

{\bf Outputs}: Estimated signal $x_{optim}$, criterion $j_t$, error $L_{optim}=\norm{y-\varphi(x_{optim})}^2$.

\section{The Cramer-Rao Lower Bounds\label{sec4}}

Consider the noisy measurement model (\ref{eq:model}), $y=\varphi(x) + \nu$, with $\nu\sim{\cal N}(0,\sigma^2)$.
Fix a direction in $H$, say $e$. We make the following two assumptions regarding the unknown signal $x$:
(1) We assume $x$ is not orthogonal to $e$, that is $\ip{x}{e}\neq 0$; (2) We assume we are given the sign of this scalar
product; say $\ip{x}{e}>0$. These two assumptions allow us to resolve the global sign ambiguity. Thus $x\in S\subset H$
where $S$ is a half-space of $H$. Since it is a convex cone we can compute expectations of random variables defined in $S$.
The likelihood function is given by
\[ L(x) = p(y|x) = \frac{1}{(2\pi)^{m/2}\sigma^m}e^{-\frac{1}{2\sigma^2}\norm{y-\varphi(x)}^2} \]
The Fisher information matrix $\I(x)$ is given by
\[ (\I(x))_{k,j} = -\E \left[
\frac{\partial^2 log\, L(x)}{\partial x_k \partial x_j}
\right]
\]
where the expectation is over the noise process, for fixed $x$.

In the following we perform the computations in the real case $H=\R^n$. For ease of notation we
assume the canonical basis of $\R^n$ and the lower index (or second index) denotes the coordinate with respect
to this basis; for instance $x_k$ and $f_{l,k}$ denote the $k^{th}$ coordinate of $x$ and $f_l$, respectively.
\[ \frac{\partial(-log\,L(x))}{\partial x_k} = \frac{1}{\sigma^2}\ip{\varphi(x)-y}{\frac{\partial \varphi (x)}{\partial x_k}} =
\frac{2}{\sigma^2} \sum_{l=1}^m \ip{x}{f_l}\left(|\ip{x}{f_l}|^2-y_l\right)f_{l,k} \]
\[ \frac{\partial^2 (-log\,L(x))}{\partial x_k\partial x_j} = \frac{2}{\sigma^2} \sum_{l=1}^m
f_{l.j} \left(|\ip{x}{f_l}|^2-y_l \right) f_{l,k} + \frac{4}{\sigma^2}\sum_{l=1}^m |\ip{x}{f_l}|^2 f_{l,j}f_{l,k} \]
Now use $\E[y_l] = |\ip{x}{y_l}|^2$. We thus obtain
\begin{equation}
\label{eq:Fisher}
\I(x) = \frac{4}{\sigma^2} \sum_{l=1}^m |\ip{x}{f_l}|^2 f_l f_l^T = \frac{4}{\sigma^2} R(x)
\end{equation}
where $R(x)$ denotes the quadratic form introduced in (\ref{eq:q2}).
Now we are ready to state the first lower bound result (see e.g. \cite{Kay2010} Theorem 3.2).

\begin{thm} \label{th4.1} In the real case $H=\R^n$, the Fisher information matrix for model (\ref{eq:model})
is given by $\I(x)$ in (\ref{eq:Fisher}). Consequently the covariance matrix of any uniabsed estimator $\Phi(y)$ for $x$
is bounded below by the Cramer-Rao lower bound as follows
\begin{equation}\label{eq:CRLB}
Var[\Phi(y)] \geq CRLB(x) := (\I(x))^{-1} = \frac{\sigma^2}{4} (R(x))^{-1}
\end{equation}
Furthermore the conditional mean-square error of any unbiased estimator $\Phi(y)$ is given by
\begin{equation}
\label{eq:mse}
\E[\norm{\Phi(y)-x}^2|x] \geq \frac{\sigma^2}{4} trace\{(R(x))^{-1}\}.
\end{equation}
When signal $x$ is random and drawn from $x\sim {\cal N}(0,I)$, the mean-square error of the unbiased
estimator $\Phi(y)$ is bounded below by
\begin{equation}
\label{eq:msse}
\E[\norm{\Phi(y)-x}^2] \geq \frac{\sigma^2}{4} trace\{ \E[(R(x))^{-1}] \}.
\end{equation}
\end{thm}
 \begin{remark}
Corollary \ref{th2.5} implies that when $\varphi$ is injective the Fisher information matrix is invertible, hence a bounded CRLB, and the
signal is identifiable in $S$ (up to a global phase factor). 
\end{remark}

We derive now a different lower bound for a modified estimation problem. Let us denote by $X=xx^T$
and $F_k=f_kf_k^T$ the rank-1 operators associated to vectors $x$ and $f_k$ respectively. Note $|\ip{x}{f_k}|^2 = trace\{XF_k\}$.
Hence
\[ y_k=trace\{XF_k\} + \nu_k~~,~~1\leq k\leq m \]
We would like to obtain a lower bound on conditional mean-square error of an unbiased estimator of the rank-1 matrix $X$.
A naive computation of the Fisher information associated to $X$ in the linear model above would produce a singular matrix whenever
$m <\frac{n(n+1)}{2}$ (the reason being the fact that a general symmetric $X$ is not identifiable merely from $m<\frac{n(n+1)}{2}$
measurements). Instead the bound should be derived under the additional hypothesis that $X$ has rank one. We obtain such a bound using 
a modified CRLB. Let $g:S\rightarrow\R^{n^2}$ be the vector valued map 
\[ g(x) = \left( \begin{array}{c}
\mbox{$x_ix_j$} 
\end{array} \right) _{1\leq i<j\leq n} \]
of $n^2$ components.
Let $\Psi(y)$ denote any unbiased estimator of the rank-1 matrix $X$. Then (see equation (3.30) in \cite{Kay2010})
\begin{equation}
\label{eq:m4}
{\cal Cov}(\Psi(y)) \geq \left(\frac{\partial g}{\partial x}\right) \I^{-1} \left(\frac{\partial g}{\partial x}\right)^T
\end{equation}
Taking trace on both sides we get
\[ \E[\norm{\Psi(y)-xx^T}^2|x] \geq trace\{\I^{-1}  \left(\frac{\partial g}{\partial x}\right)^T \left(\frac{\partial g}{\partial x}\right) \}
\]
Let $H=\frac{\partial g}{\partial x}$. Then we have
\begin{eqnarray*}
 (H^TH)_{k_1,k_2} & = & \sum_{j}H_{j,l_1}H_{j,l_2} = \sum_j \frac{\partial g_j}{\partial x_{l_1}}\frac{\partial g_j}{\partial x_{l_2}} \\
& = & \sum_{i,j=1}^n \frac{\partial (x_ix_j)}{\partial x_{l_1}}\frac{\partial (x_ix_j)}{\partial x_{l_2}} = 2\norm{x}^2 \delta_{l_1,l_2}+
2x_{l_1}x_{l_2}
\end{eqnarray*}
Thus we obtained
\begin{thm}\label{th4.1b}
The conditional mean-square error of any unbiased estimator of the rank-1 matrix $X=xx^T$ is bounded below by
\begin{equation}\label{eq:M4}
\E[\norm{\Psi(y)-xx^T}^2|x] \geq \frac{\sigma^2}{2}\left( \norm{x}^2 trace\{R^{-1}\} + x^T R^{-1}x \right). 
\end{equation}
\end{thm}

Consider now the case of the Maximum Likelihood Estimator (MLE) whose optimization problem was considered in
the previous section. For model (\ref{eq:model}) this takes the form of
\begin{equation}
\label{eq:MLE}
 \Phi_{MLE}(y) = argmin_u \norm{\varphi(u)-y}^2 = argmin_u J(u,u;0,0) 
\end{equation}
The MLE computes the global minimum in the optimization problem above. Assume that 
$\Phi_{MLE}$ selects the closest global minimum to $x$. We want to estimate lower bounds on the MLE performance
so that we can benchmark performance of any optimization algorithm against these bounds.

First we estimate the bias of the MLE estimator in the asymptotic limit
$\sigma\rightarrow 0$. The estimator must satify the MLE equation
\[ \nabla (log\,L(x)){|}_{x=\Phi_{MLE}(y)} = 0 \]
which turns into
\[ \sum_{k=1}^m \ip{\Phi_{MLE}(y)}{f_k}\left( |\ip{\Phi_{MLE}(y)}{f_k}|^2-y_k \right) f_k = 0 \]
Denote $\deltavect = \Phi_{MLE}(y)-x$. The bias is given by $Bias(x) = \E[\deltavect|x]$. Assymptotically we can assume
$\norm{\deltavect}$ is small with high probability. We shall expand the products in the above equation
taking into account only the first terms in $\deltavect$:
\[ \sum_{k=1}^m (\ip{x}{f_k} + \ip{\deltavect}{f_k})(|\ip{x}{f_k}+\ip{\deltavect}{f_k}|^2-y_k)f_k = 0\]
Expanding the products and neglecting higher order terms in $\deltavect$ we obtain:
\begin{equation}
\label{eq:v}
 \left( \sum_{k=1}^m (3 |\ip{x}{f_k}|^2-y_k) f_k f_k^T \right) \deltavect + \sum_{k=1}^m (|\ip{x}{f_k}|^2-y_k)\ip{x}{f_k}f_k = 0 
\end{equation}
Note $y_k=|\ip{x}{f_k}|^2+\nu_k$. Let us denote
\[ N = \sum_{k=1}^m \nu_k f_kf_k^T ~~,~~ R=R(x)=\sum_{k=1}^m |\ip{x}{f_k}|^2 f_kf_k^T \]
The equation that $v$ satisfies becomes $(2R-N)v-Nx=0$. Therefore
\[ v = (2R-N)^{-1}Nx~~\Rightarrow~~Bias(x) = \E[ (2R-N)^{-1}N]x. \]
For fixed $x\neq 0$, due to the lower bound in (\ref{eq:q2}) we obtain with high probability $\norm{N} < 2 eigmin(R)$,
 where $eigmin(R)$ denotes the smallest eigenvalue of $R$. Note $eigmin(R)=\frac{1}{\norm{R^{-1}}}\geq a_0\norm{x}^2$ by (\ref{eq:a0}).
Then using Neumann's series expansion we get
\[ \left( \sum_{k=1}^m (3|\ip{x}{f_k}|^2-y_k)f_kf_k^T \right)^{-1} = (2R-N)^{-1} = 
\frac{1}{2}R^{-1}+\frac{1}{4}R^{-1}NR^{-1} + O(N^2) \]
Thus we obtain
\[ \deltavect = (2R-N)^{-1} Nx = \frac{1}{2}R^{-1}Nx + \frac{1}{4}R^{-1}NR^{-1}Nx + O(N^3) \]
Note also the similarity criterion in expansion above is related to
$\frac{\norm{N}}{eigmin(R)}$ which is of the order $\sigma\norm{R^{-1}}$.
Since all odd moments of Gaussian random variables vanish we obtain $\E[(R^{-1}N)^{2p+1}] = 0$. Hence
\begin{eqnarray}
 Bias(x) & = & \E[\deltavect] = \frac{1}{4} \E[R^{-1}NR^{-1}Nx] +\E[O(N^4)] =  \nonumber \\
  &  = &  \frac{\sigma^2}{4}\sum_{k=1}^m \ip{x}{f_k} \ip{R^{-1}f_k}{f_k}R^{-1}f_k + O((\sigma\norm{R^{-1}})^4) \label{eq:bias}
\end{eqnarray}
The leading term in bias has the form
\begin{equation} 
\label{eq:bias0}
Bias^0(x) = \frac{\sigma^2}{4}\delta~~,~~{\rm where}~~\delta=\sum_{k=1}^m \ip{x}{f_k}\ip{(R(x))^{-1}f_k}{f_k}(R(x))^{-1}f_k 
\end{equation}
Note the dependence on $x$ is highly nonlinear. 
We would like next to obtain the modified CRLB for MLE taking into account its bias. We need to estimate the first 
derivatives of $Bias(x)$ with respect to $x$, $\frac{\partial (Bias(x))_j}{\partial x_l}$.
Again we shall derive the asymptotic approximation of this matrix: 
\[ \Delta(x)_{j,l} = \frac{4}{\sigma^2} \frac{\partial (Bias^0(x))_j}{\partial x_l} = \frac{\partial \delta_j(x)}{\partial x_l}. \]
The key relation to use is
\[ \frac{\partial}{\partial x_l}  (R(x))^{-1}f_k = - R^{-1}\frac{\partial R}{\partial x_l}R^{-1}f_k=-2
\sum_{p=1}^m \ip{x}{f_p}f_{p,l} \ip{R^{-1}f_k}{f_p}R^{-1}f_p \]
which comes from $R(x)(R(x))^{-1}f_k = f_k$ by differentiating with respect to $x_l$, and from (\ref{eq:q2}).
After some straightforward but tedious algebra we obtain

\begin{eqnarray}
\label{eq:Delta}
\Delta & = &  \sum_{k=1}^m \ip{R^{-1}f_k}{f_k}R^{-1}f_kf_k^T   \\
 & - & 2\sum_{k,p=1}^m \ip{x}{f_k}\ip{x}{f_p}\ip{R^{-1}f_k}{f_k}\ip{R^{-1}f_k}{f_p}R^{-1}f_pf_p^T \nonumber \\
 & - & 2\sum_{k,p=1}^m \ip{x}{f_k} \ip{x}{f_p}| \ip{R^{-1}f_p}{f_k} |^2 R^{-1}f_k f_p^T  \nonumber
\end{eqnarray}

Now we can compute the modified Cramer-Rao lower bound for the MLE estimator (see e.g. \cite{Kay2010} Equation (3.30)). 
\begin{thm}\label{th4.2} The MLE estimator (\ref{eq:MLE}) is biased. Its expectation admits the following asymptotic approximation
\begin{equation}
\label{eq:meanMLE}
\E[\Phi_{MLE}(y)] = x + \frac{\sigma^2}{4}\delta + O((\sigma\norm{R^{-1}})^4).
\end{equation}
Its covariance matrix is bounded below by
\begin{equation}
\label{eq:covMLE}
{\cal Cov}[\Phi_{MLE}(y)] \geq (I+\frac{\partial Bias}{\partial x})\I^{-1}(I+(\frac{\partial Bias}{\partial x})^T) = \frac{\sigma^2}{4}R^{-1} + \frac{\sigma^4}{16}(R^{-1}\Delta^T+\Delta R^{-1})
+O((\sigma\norm{R^{-1}})^6)
\end{equation}
where $I$ is the identity matrix. Furthermore, the conditional mean-square error is bounded below by
\begin{eqnarray}
\label{eq:mseMLE} \\
\E[\norm{\Phi_{MLE}(y)-x}^2|x] & =& \norm{Bias(x)}^2 + trace\{{\cal Cov}[\Phi_{MLE}(y)]\} \nonumber \\
\geq \frac{\sigma^2}{4}trace\{R^{-1}\} & + &
\frac{\sigma^4}{16}\left( \norm{\delta}^2+2trace\{\Delta R^{-1} \} \right) + O((\sigma\norm{R^{-1}})^6) \nonumber
\end{eqnarray}
Here we used the notation $R=R(x)=\sum_{k=1}^m|\ip{x}{f_k}|^2 f_kf_k^T$, and $\delta,\Delta$ given by (\ref{eq:bias0}), (\ref{eq:Delta}).
\end{thm}

\section{Numerical Analysis\label{sec5}}

In this section we present numerical simulations for the algorithms presented in this paper.

We generated random frames or redundancy 3, that is $m=3n$, as well as random signals $x$.
All these vectors (frame and signal) are drawn from ${\cal N}(0,I)$. We set the first component of $x$
positive, and so we decided the global sign after reconstruction. 
To the magnitude square of signal coefficients $\varphi(x)$ we added Gaussian noise of variance $\sigma^2$ 
to achieve a fixed signal-to-noise-ratio defined as
\[ SNR = \frac{\sum_{k=1}^m|\ip{x}{f_k}|^4}{m\sigma^2} ~~,~~SNRdB = 10\log_{10} (SNR)~~{\rm [dB]}\]
Note the similarity criterion used in asymptotic expansions (\ref{eq:meanMLE}) and (\ref{eq:covMLE}) is of the same
order as $\sigma\norm{R^{-1}} \sim \frac{1}{\sqrt{SNR}}$ (up to multiplicative constants).
We used 11 values of SNRdB in 10dB increments from -20dB to +80dB. 

For the first algorithm, results are averaged over 100 noise realizations.
In each instance of the algorithm we initialized $x^0,\lambda_0,\mu_0$ as described in subsection \ref{subsec3.4}
with $\alpha=0.9$. At each iteration $\lambda_{t+1}=\lambda_t/1.05$. We run the algorithm for at least 100 steps, or
until $\lambda_t$ gets below $10^{-8}$. The parameter$\mu_t$ is adapted as follows:
 $\mu_t=max(1,\lambda_t)$.

\begin{figure}[htb]
\includegraphics[width=100mm,height=60mm]{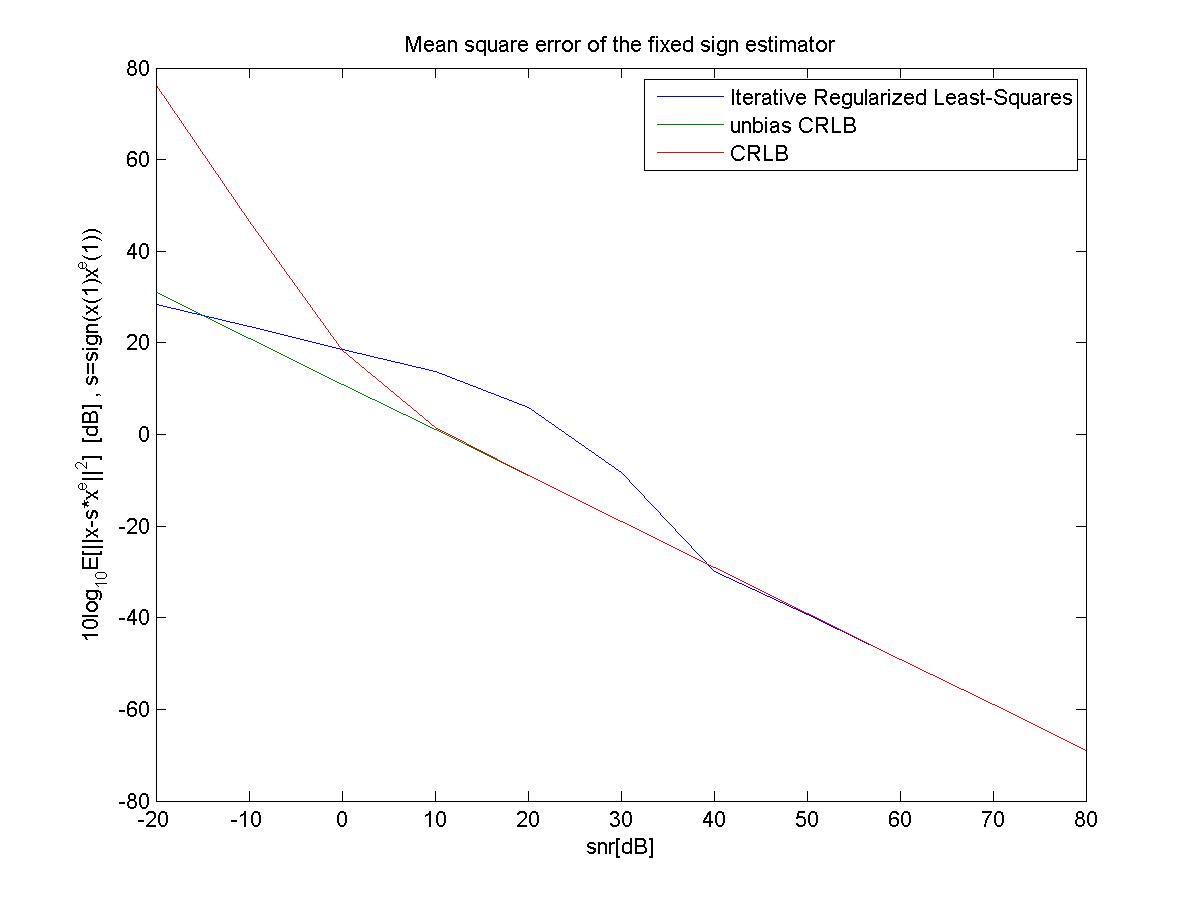}
\includegraphics[width=100mm,height=60mm]{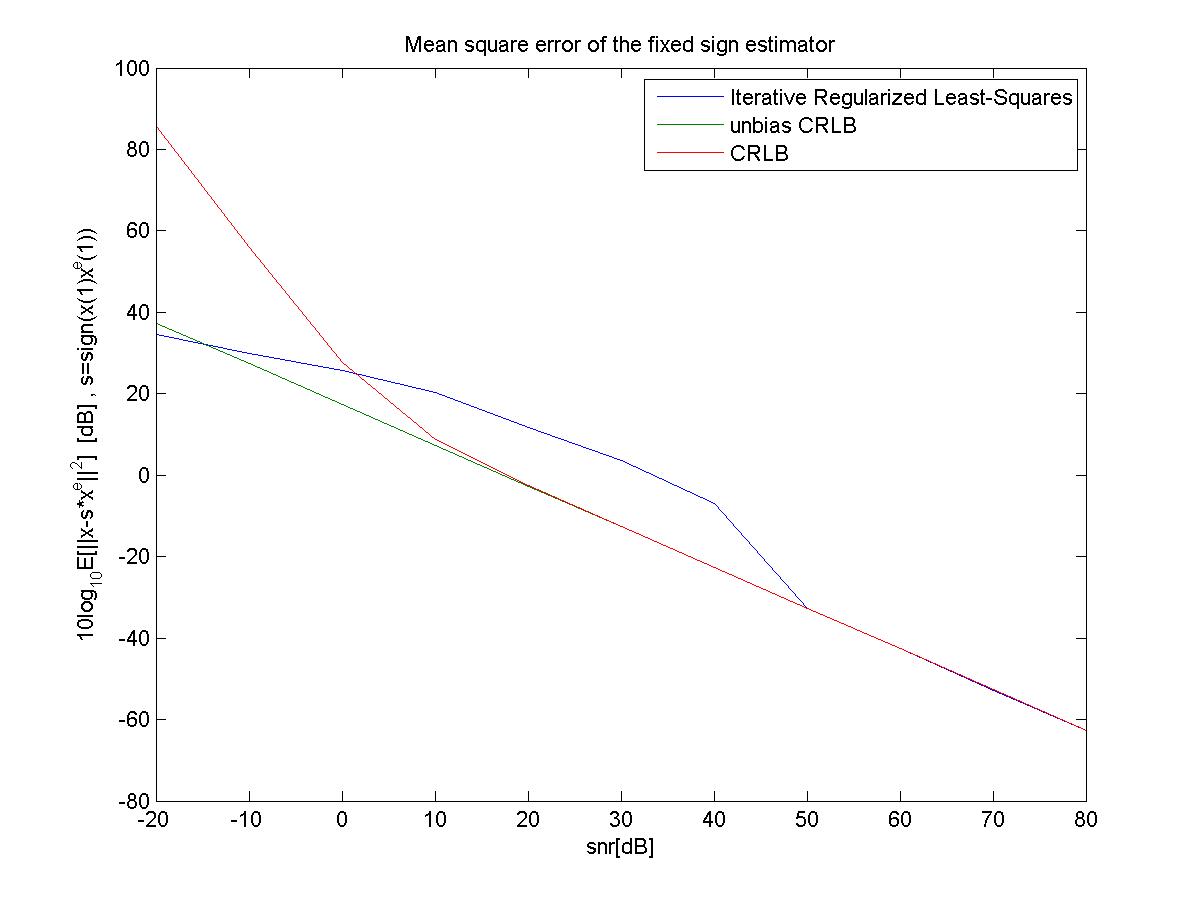}
\includegraphics[width=100mm,height=60mm]{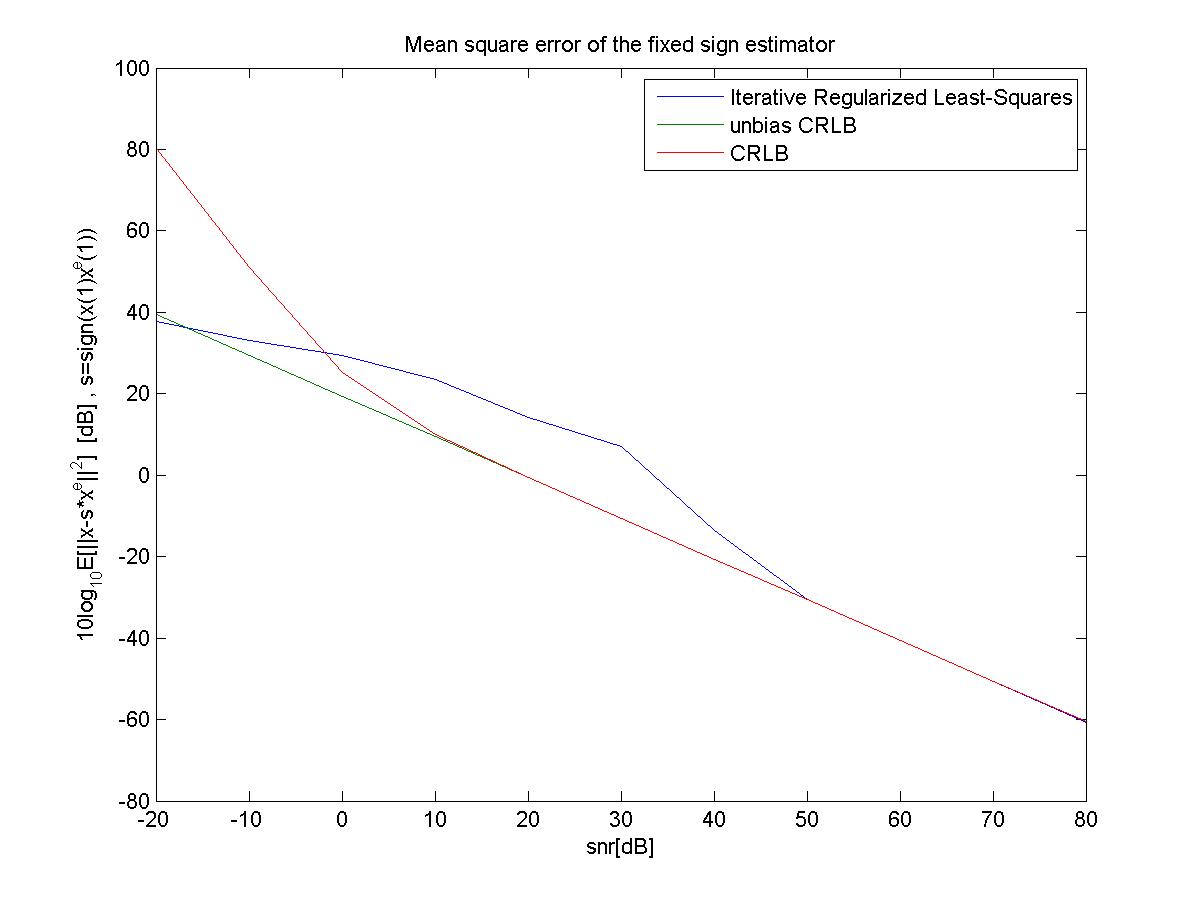}
\caption{Conditional Mean-Square Error and CRLB bounds for $n=10$ (top plot), $n=50$ (middle plot),
 and $n=100$ (bottom plot). \label{fig1}}
\end{figure}

We include results for $n=10$, $n=50$ and $n=100$.
Figure \ref{fig1} includes the conditional mean-square error averaged over 100 noise realizations, and the lower
bounds: the unbiased CRLB (\ref{eq:CRLB}) and the MLE adapted CRLB (\ref{eq:mseMLE}). 
Note the two lower bounds are indistingueshable for $SNR>20dB$. For low SNR, when the
two bounds differ significantly, the approximation (\ref{eq:v}) is no longer valid. Hence the
bound would be different as well. In general we cannot obtaion a closed form solution for the
new bound. 

\begin{figure}[htb]
\includegraphics[width=100mm,height=60mm]{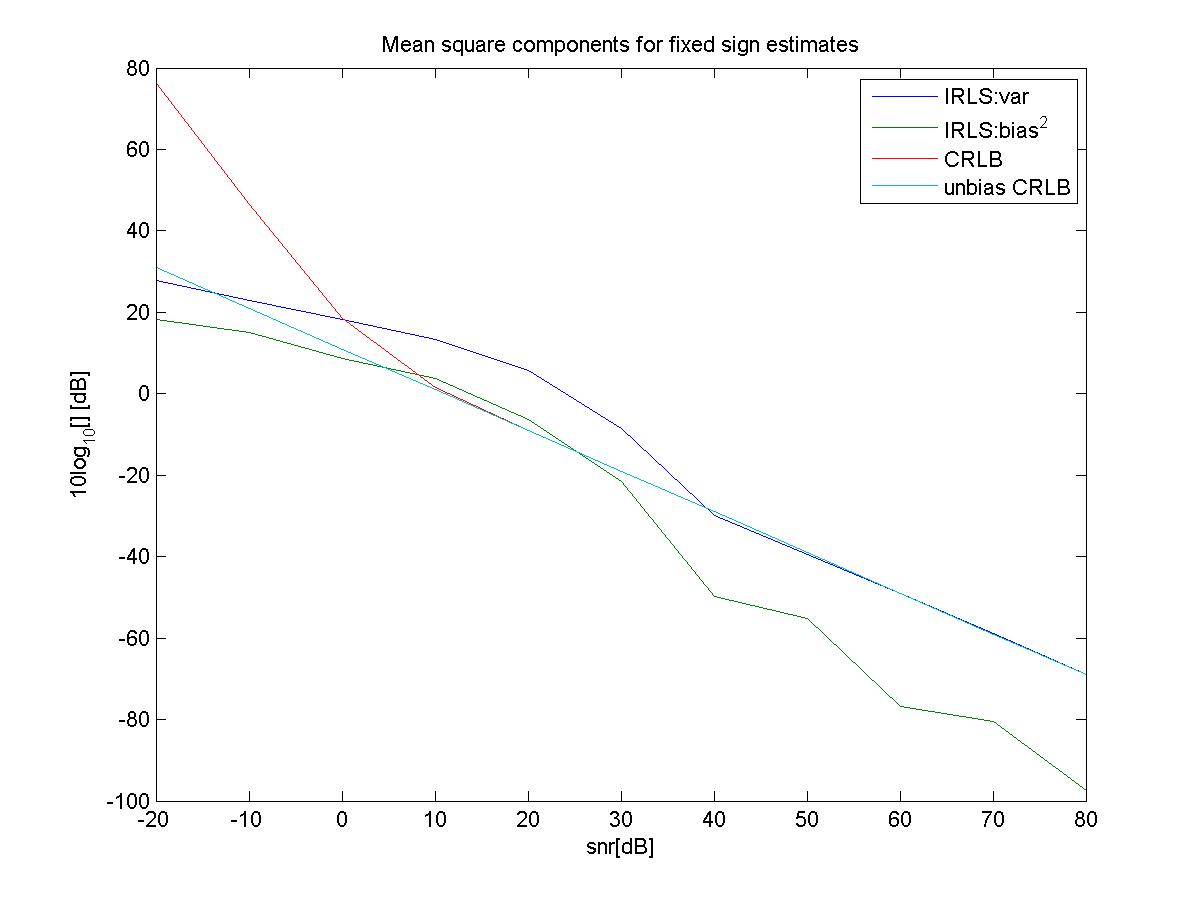}
\includegraphics[width=100mm,height=60mm]{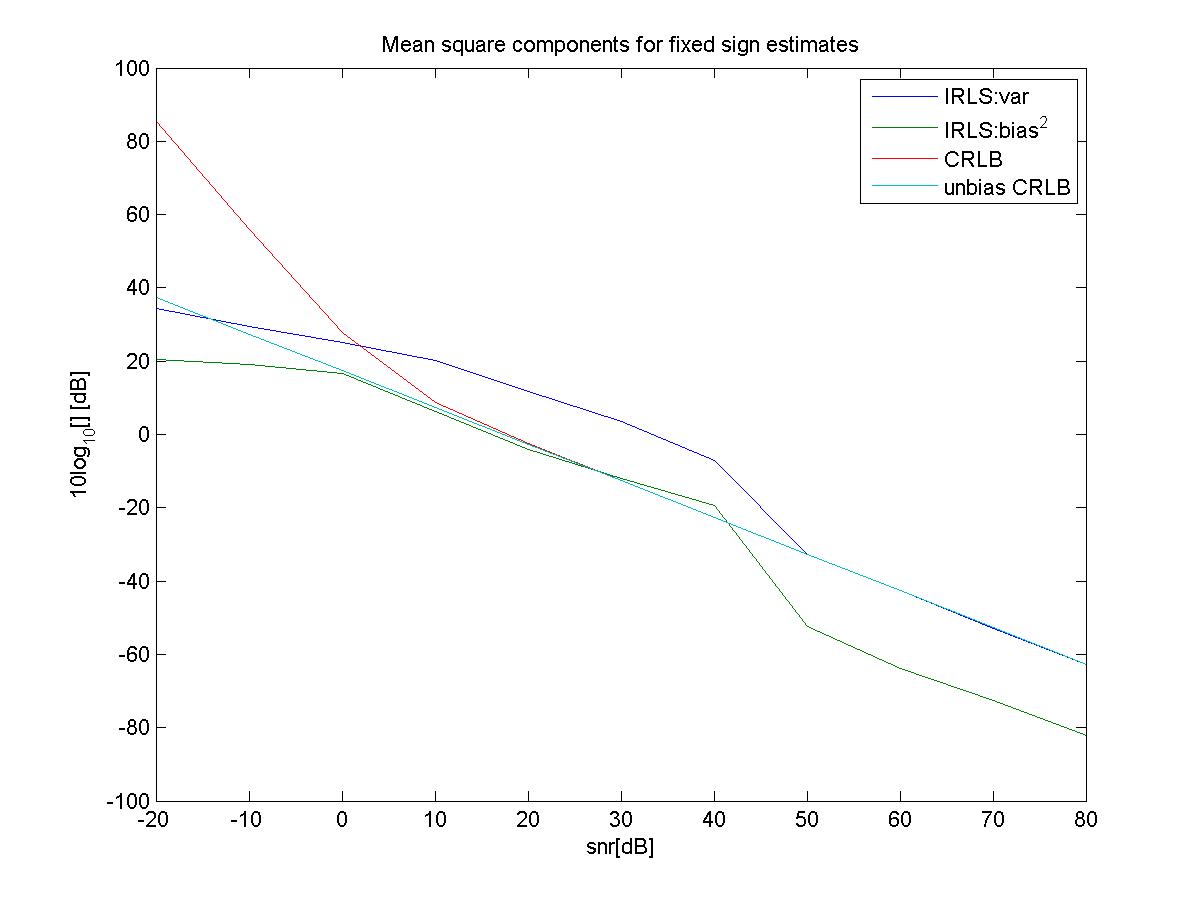}
\includegraphics[width=100mm,height=60mm]{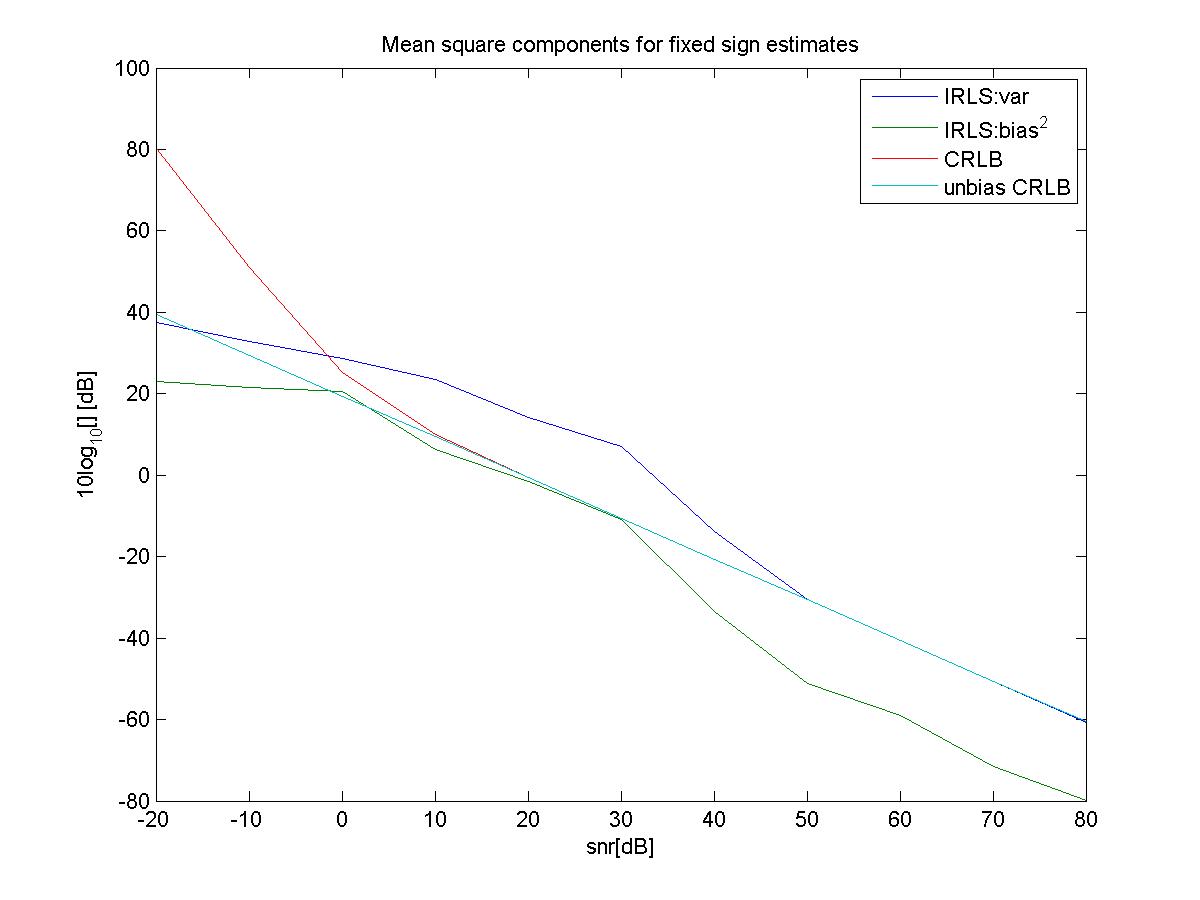}
\caption{Bias and Variance components of conditional Mean-Square Error and CRLB bounds for $n=10$ (top plot), $n=50$ (middle plot),
 and $n=100$ (bottom plot). \label{fig2}}
\end{figure}

In Figure \ref{fig2} we plot the bias and variance components of the mean-square error for the same results
in Figure \ref{fig1}. Note the bias is relatively small. The bulk of mean-square error is due to estimation variance.

\begin{figure}[htb]
\includegraphics[width=100mm,height=60mm]{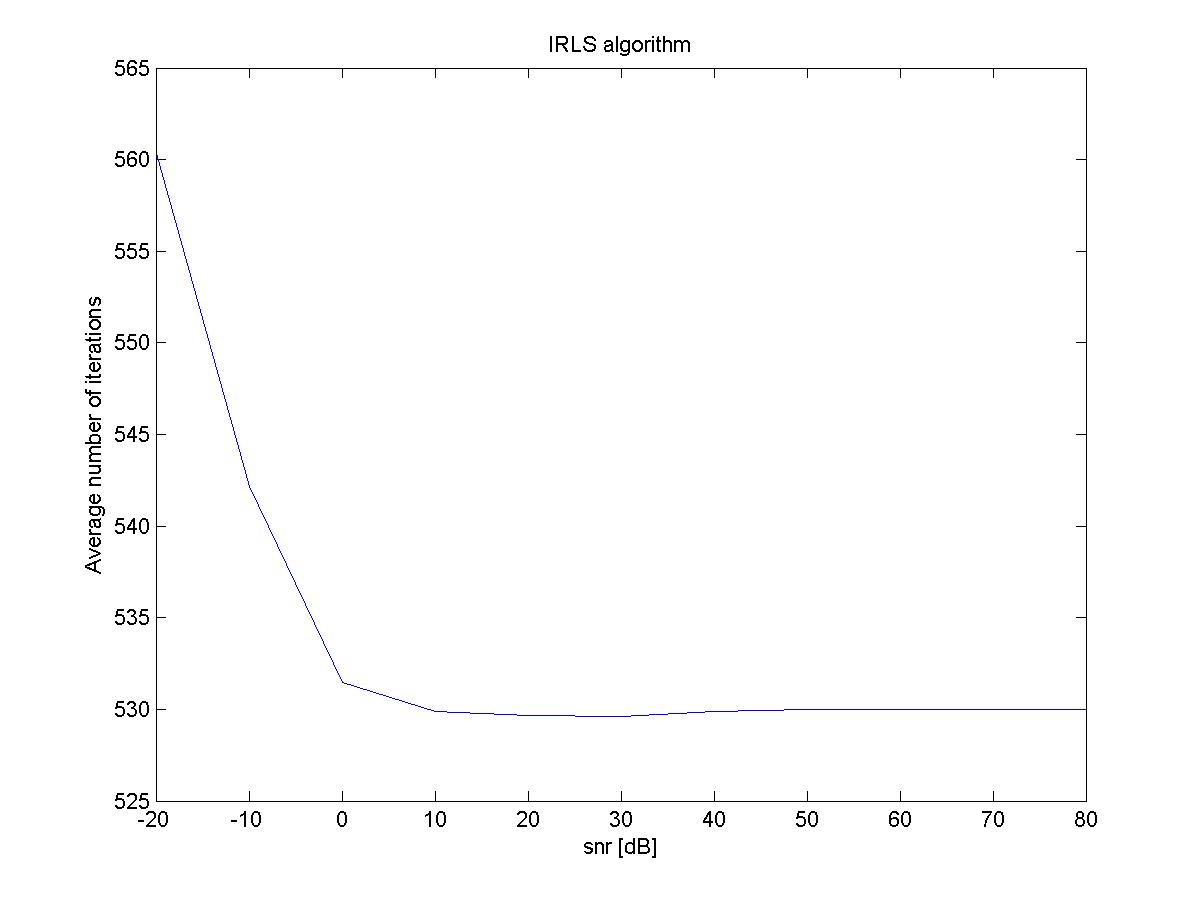}
\includegraphics[width=100mm,height=60mm]{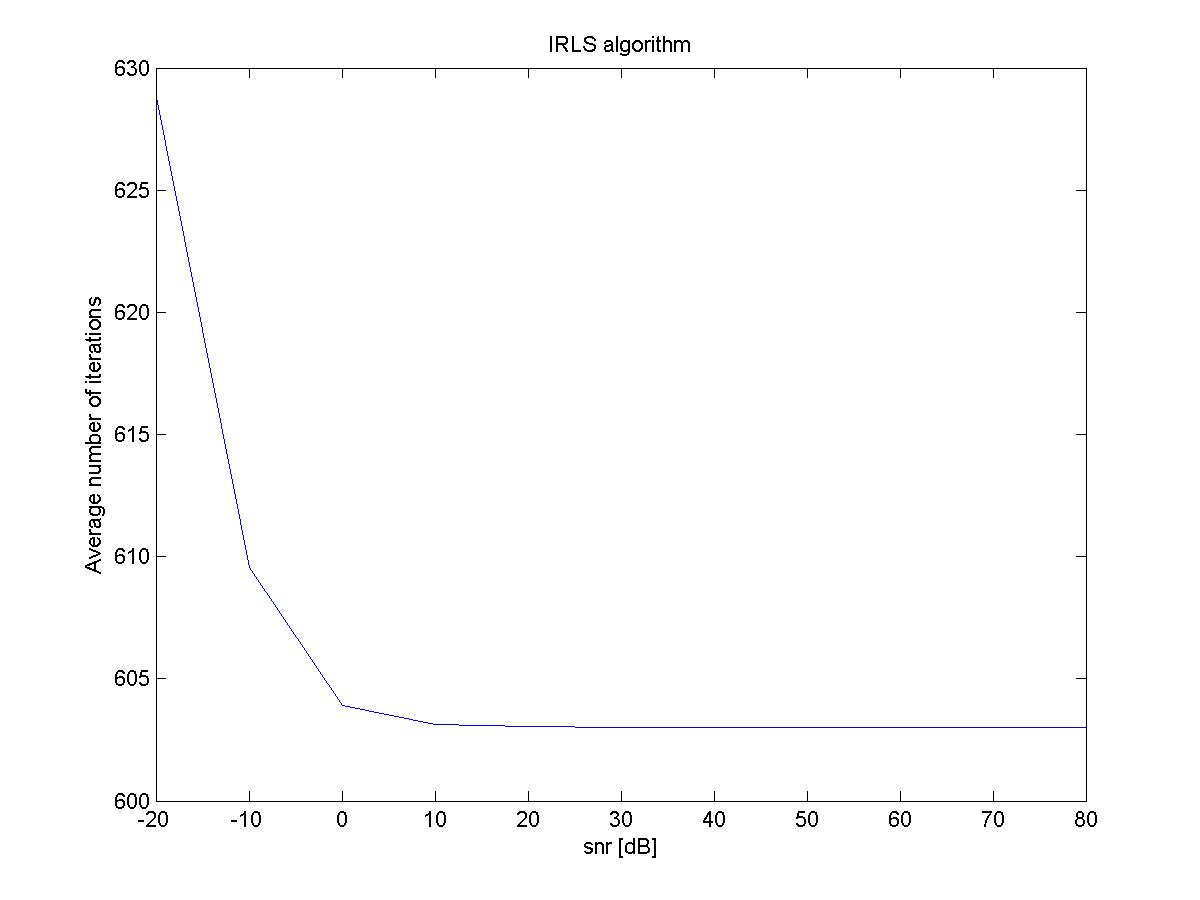}
\includegraphics[width=100mm,height=60mm]{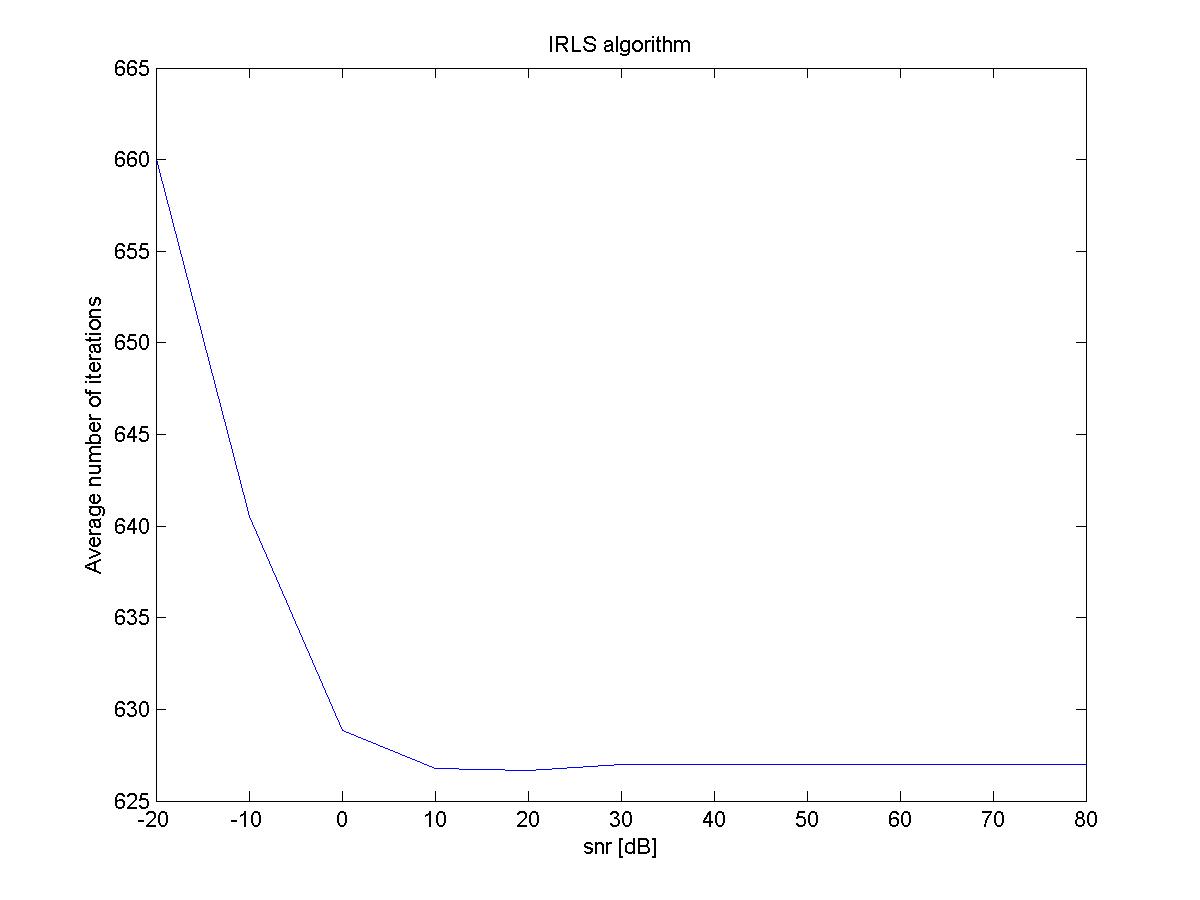}
\caption{Average number of iterations vs. SNR for $n=10$ (top plot), $n=50$ (middle plot),
 and $n=100$ (bottom plot). \label{fig3}}
\end{figure}

Figure \ref{fig3} contains the average number of iterations for each of these cases. The algorithm
runs for about 530-660 steps. 

\begin{figure}[htb]
\includegraphics[width=100mm,height=60mm]{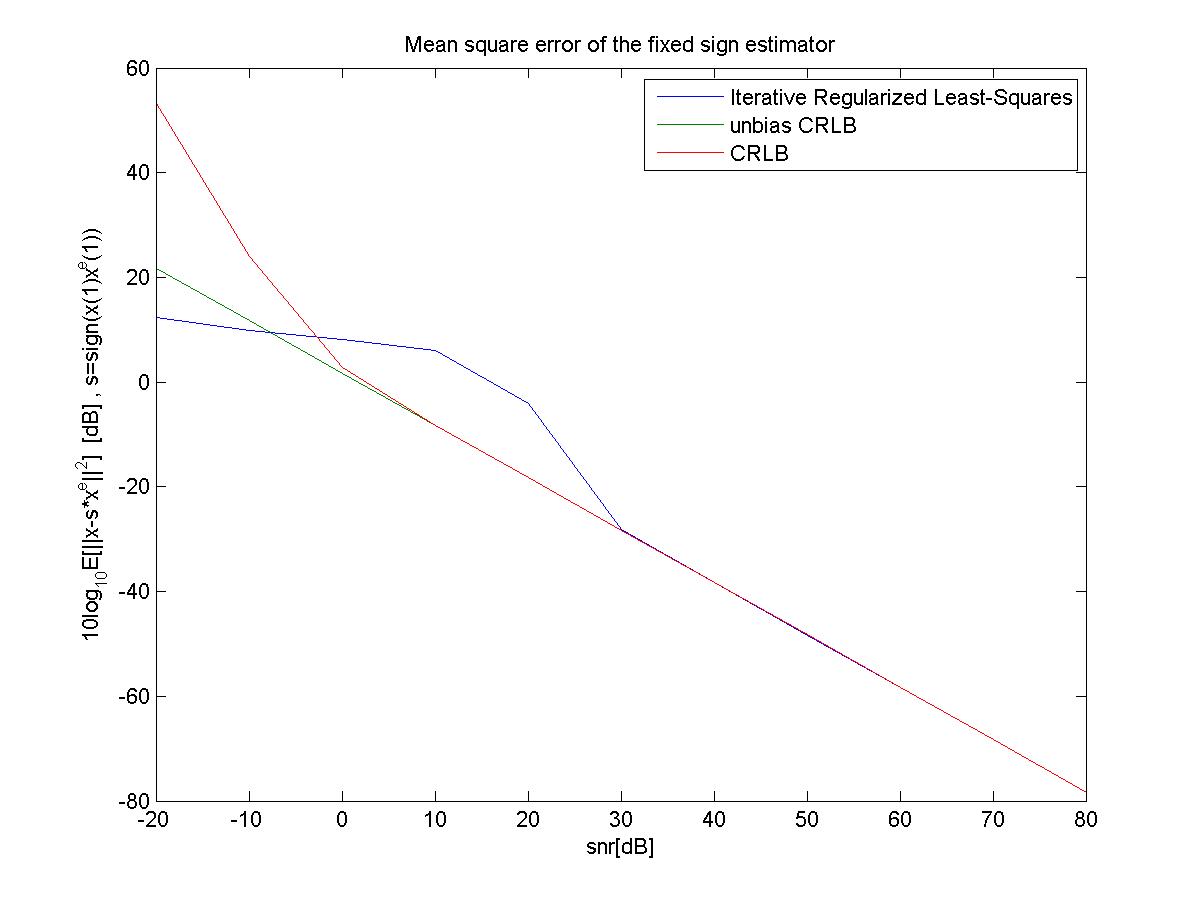}
\includegraphics[width=100mm,height=60mm]{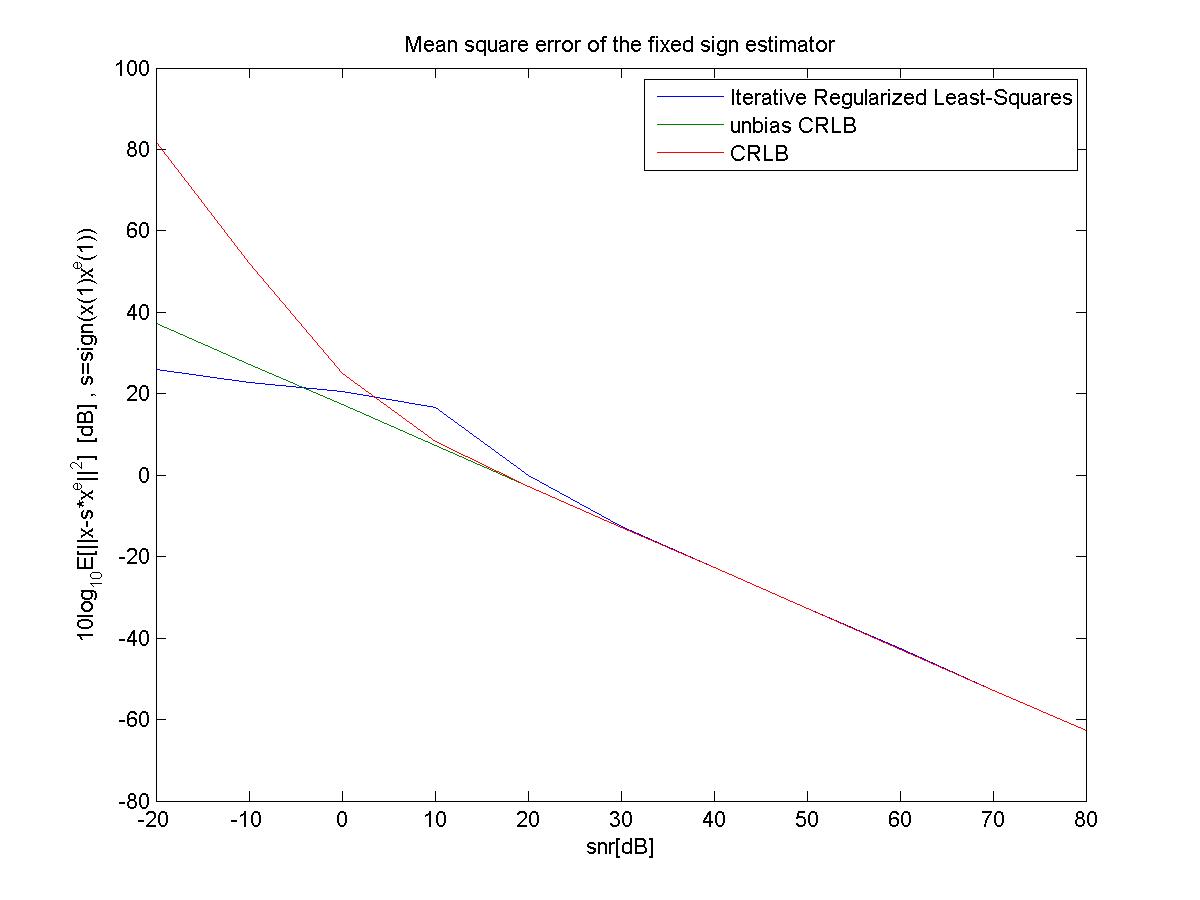}
\includegraphics[width=100mm,height=60mm]{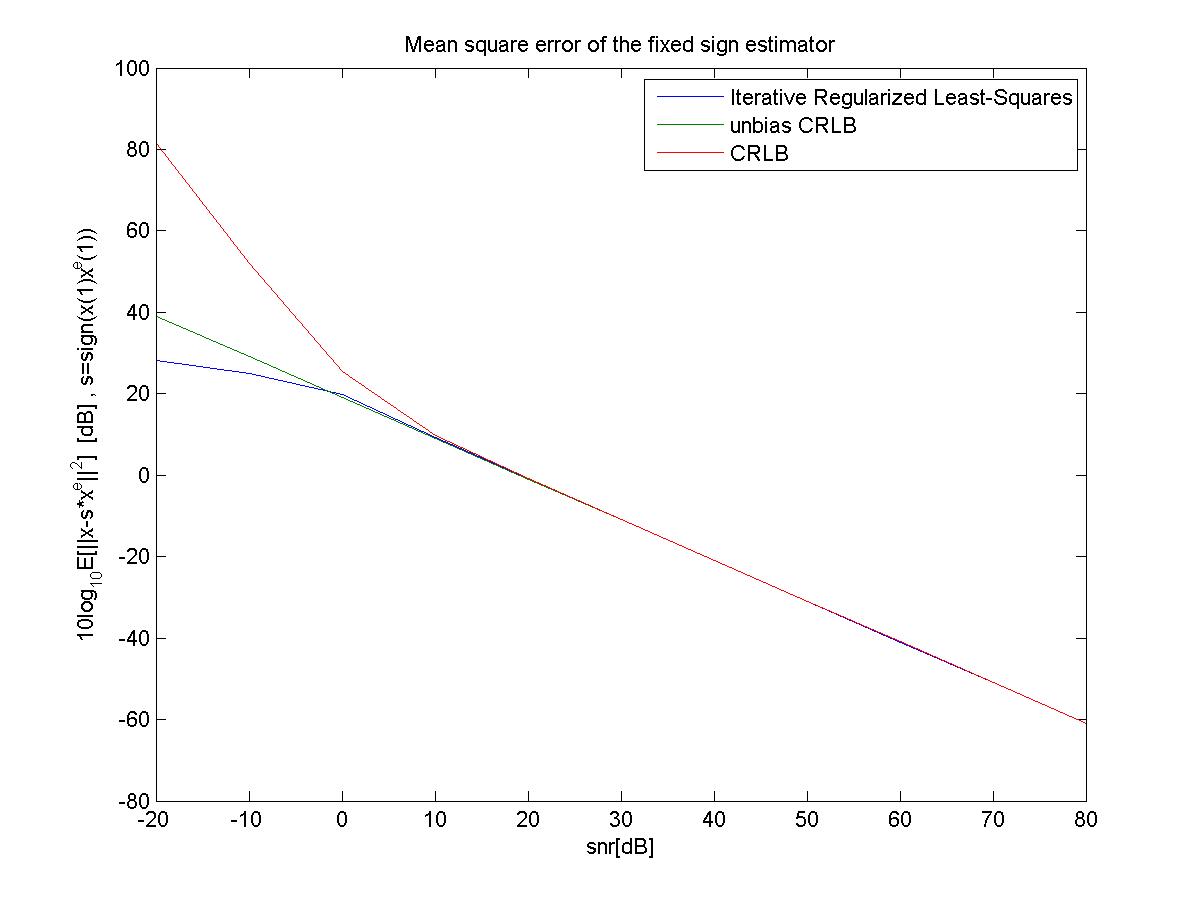}
\caption{Conditional Mean-Square Error and CRLB bounds for $n=10$ (top plot), $n=50$ (middle plot),
 and $n=100$ (bottom plot) for fixed sign. \label{fig4}}
\end{figure}

\begin{figure}[htb]
\includegraphics[width=100mm,height=60mm]{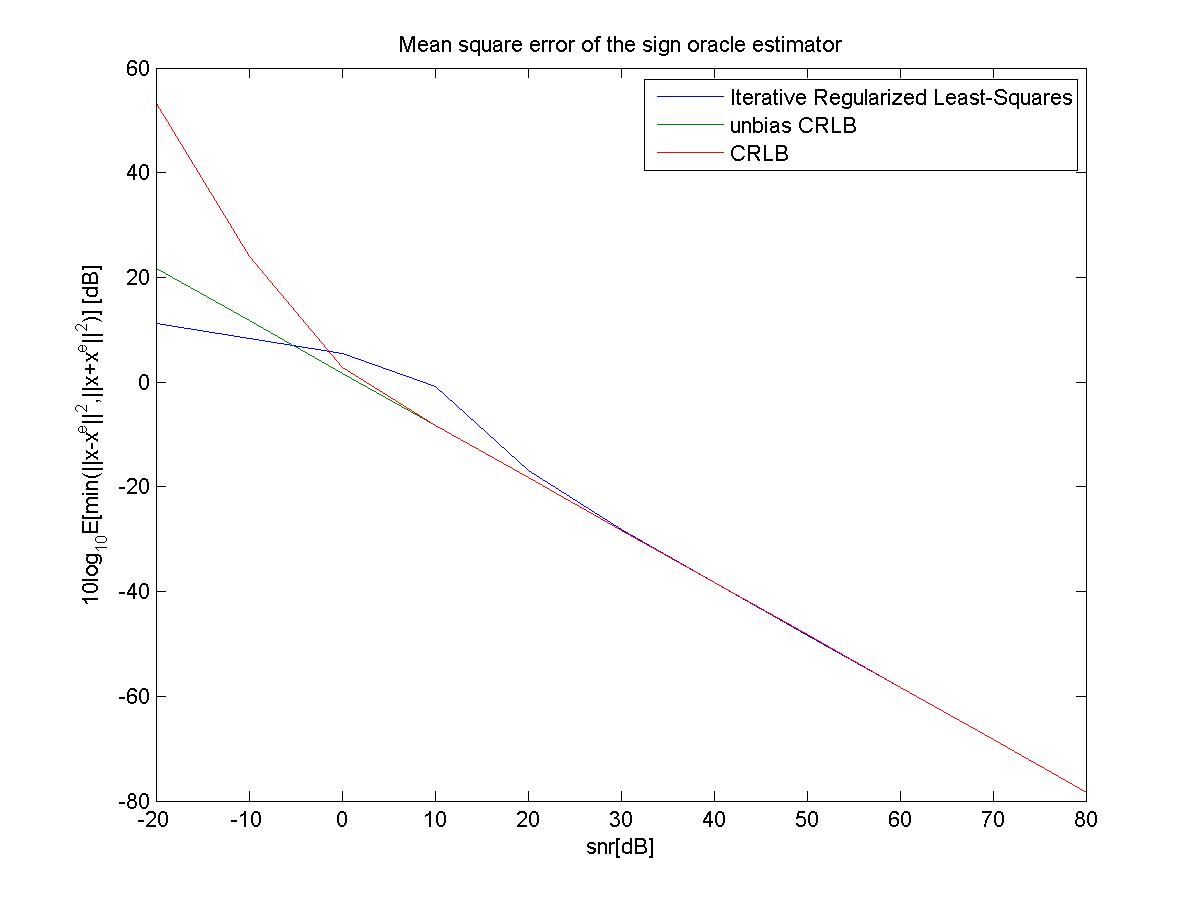}
\includegraphics[width=100mm,height=60mm]{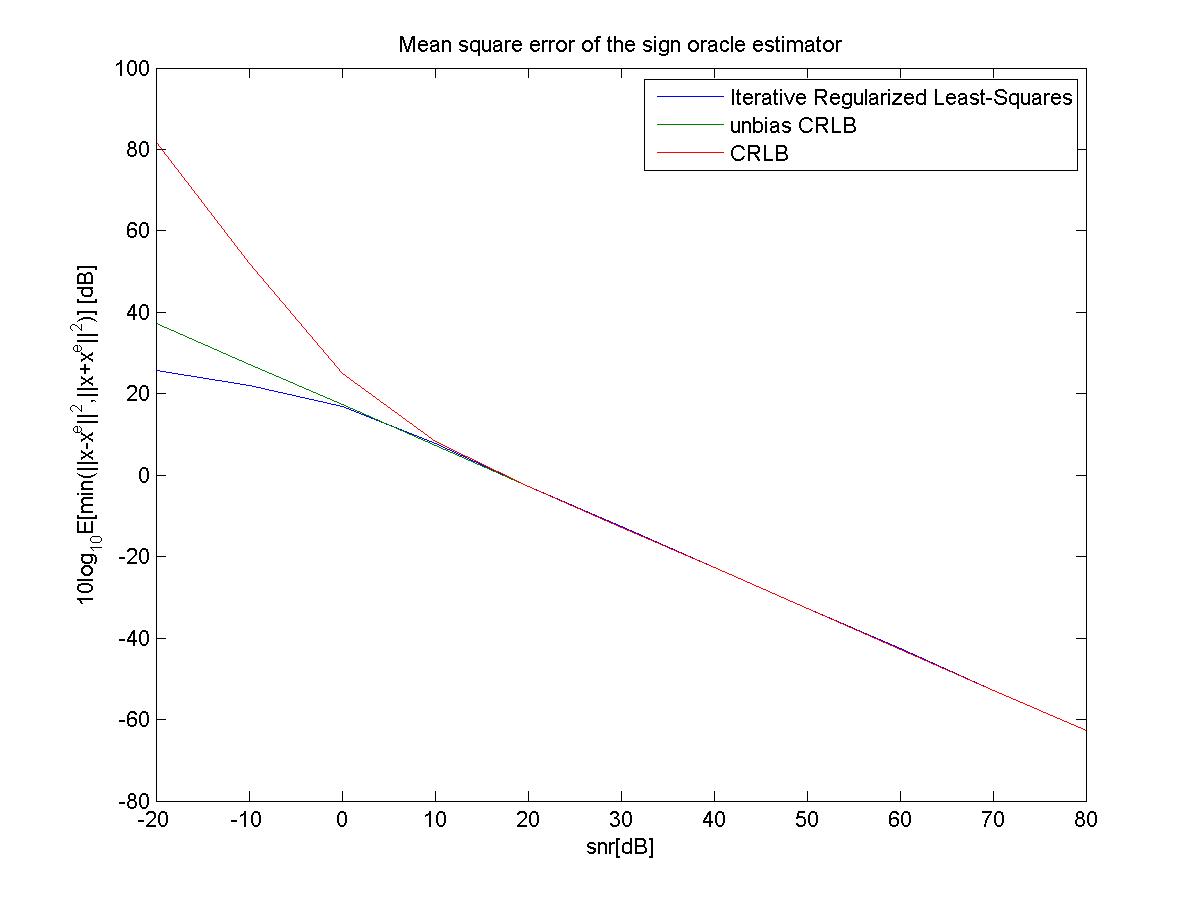}
\includegraphics[width=100mm,height=60mm]{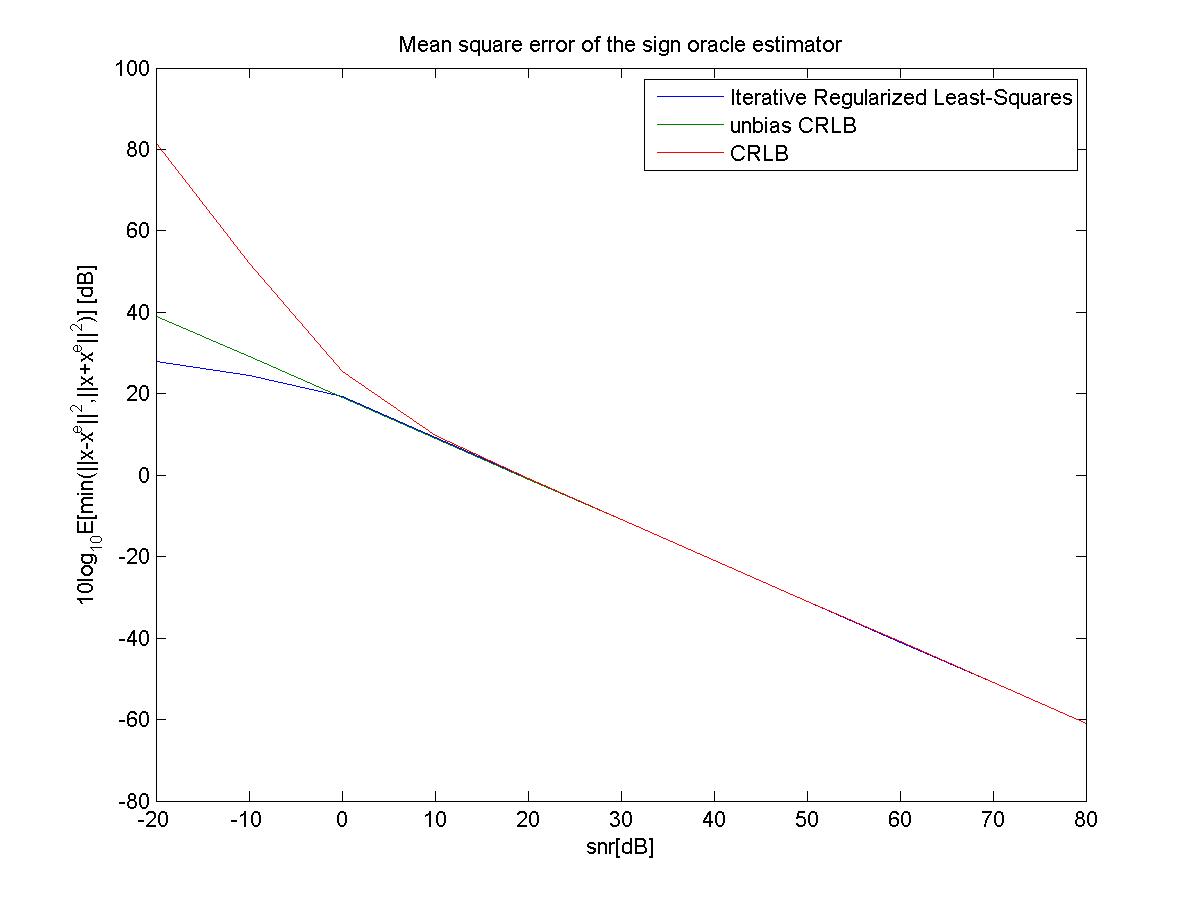}
\caption{Conditional Mean-Square Error and CRLB bounds for $n=10$ (top plot), $n=50$ (middle plot),
 and $n=100$ (bottom plot) for global sign oracle. \label{fig5}}
\end{figure}

For the second algorithm we repeated the same cases ($n=10,50,100$) and same levels of SNR, but we average over 1000
noise realizations. We present the mean-square error in two cases: in Figure \ref{fig4} the case of fixed sign as discussed
earlier (first component of $x$ is positive); in Figure \ref{fig5} the case of a sign oracle, when the global sign
is chosen as given by the minimum $min(\norm{x-x_{optim}},\norm{x+x_{optim}})$. 

\section{Conclusions\label{sec6}}
Novel necessary conditions for signal reconstruction from magnitudes of frame coefficients have been presented.
These conditions are also sufficient in the real case. The least-square criterion has been analyzed, and two
algorithms have been proposed to optimize this criterion.
Performance of the second algorithm presented in this paper is remarkably close to the
theoretical lower bound given by the Cramer-Rao inequality. In fact for low SNR its 
performance is better than the asymptotic approximation given by the modified CRLB.

\end{document}